\documentclass[12pt]{article}
\usepackage{amsmath}
\usepackage{amssymb}

\newcommand{\acr}{AC(\Rbar)}
\newcommand{\supp}{{\rm supp}}

\newcommand{\alexc}{{\cal A}_C}

\newcommand{\balexc}{{\cal B}_C}
\newcommand{\Rbar}{\overline{\R}}
\newcommand{\D}{{\cal D}}
\newcommand{\Dp}{{\cal D}'}
\newcommand{\nbv}{{\cal NBV}}
\newcommand{\ebv}{{\cal EBV}}
\newcommand{\bv}{{\cal BV}}
\newcommand{\intinf}{\int^\infty_{-\infty}}

\newcommand{\N}{{\mathbb N}}
\newcommand{\R}{{\mathbb R}}

\newcommand{\fn}{\!:\!}

\newcommand{\llim}{\lim\limits}

\newcommand{\qed}{\mbox{$\quad\blacksquare$}}
\newcounter{examples}

\newtheorem{theorem}{Theorem}

\newtheorem{prop}[theorem]{Proposition}
\newtheorem{corollary}[theorem]{Corollary}
\newtheorem{remark}[theorem]{Remark}
\newtheorem{defn}[theorem]{Definition}

\begin{document}
\hspace{-2cm}
\raisebox{12ex}[1ex]{\fbox{{\footnotesize
Preprint
September 18, 2009.  To appear in {\it Abstract and Applied Analysis}
}}}

\begin{center}
{\large\bf Convolutions with the continuous primitive integral}
\vskip.25in
Erik Talvila\footnote{Supported by the
Natural Sciences and Engineering Research Council of Canada.
}\\ [2mm]
{\footnotesize
Department of Mathematics and Statistics \\
University of the Fraser Valley\\
Abbotsford, BC Canada V2S 7M8\\
Erik.Talvila@ufv.ca}
\end{center}

{\footnotesize
\noindent
{\bf Abstract.} 
If $F$ is a continuous function on the real line and $f=F'$ is
its distributional derivative then the continuous primitive
integral of distribution $f$ is $\int_a^bf=F(b)-F(a)$.  This integral
contains the Lebesgue, Henstock--Kurzweil and wide Denjoy integrals.
Under the Alexiewicz norm the space of integrable distributions is
a Banach space.  We define the convolution $f\ast g(x)=\intinf f(x-y)g(y)\,dy$
for $f$ an integrable distribution and $g$ a function of bounded variation
or an $L^1$ function.
Usual properties of convolutions are shown  to hold: commutativity,
associativity, commutation with translation.  For $g$ of bounded
variation, $f\ast g$ is uniformly continuous and we have
the estimate $\|f\ast g\|_\infty\leq \|f\|\|g\|_\bv$ where 
$\|f\|=\sup_I|\int_If|$ 
is the
Alexiewicz norm.  This supremum is taken over all intervals  $I\subset\R$.
When $g\in L^1$ the estimate is
$\|f\ast g\|\leq \|f\|\|g\|_1$.  There are results on differentiation
and integration of  convolutions.  A type of Fubini theorem is proved
for the continuous primitive integral.
\\
{\bf 2000 subject classification:} 26A39, 42A85, 46E15, 46F10, 46G12
}

\section{Introduction and notation}\label{introduction}
The convolution of two functions $f$ and $g$ on the real line is
$f\ast g(x)=\intinf f(x-y)g(y)\,dy$.  
Convolutions play
an important role in pure and applied mathematics in Fourier analysis,
approximation theory, differential equations, integral equations
and many other areas. 
In this paper we consider convolutions
for the continuous primitive integral.
This integral extends  the Lebesgue, Henstock--Kurzweil
and wide Denjoy integrals on the real line and has a very simple 
definition in terms
of distributional derivatives.

Some of the main results for Lebesgue integral convolutions are that
the convolution defines a Banach algebra  on $L^1$ and
$\ast\fn L^1\times L^1\to L^1$ such that $\|f\ast g\|_1\leq \|f\|_1 \|g\|_1$.
The convolution is commutative, associative and
commutes with translations. If $f\in L^1$ and $g\in C^n$ then
$f\ast g\in C^n$ and $(f\ast g)^{(n)}(x)=f\ast g^{(n)}(x)$.  Convolutions
also have the approximation property that if $f\in L^p$ ($1\leq p<\infty$)
and $g\in L^1$ then $\|f\ast g_t-af\|_p\to 0$ as $t\to 0$, where
$g_t(x)=g(x/t)/t$ and $a=\intinf g$.  
When $f$ is bounded and continuous, there is a 
similar result for $p=\infty$.  For these results see, for example, 
\cite{folland}.
See \cite{talvilafourier} for related results with the Henstock--Kurzweil
integral.
Using the Alexiewicz norm, all of these results have generalisations to
continuous primitive integrals that are proven below.

We now define the  continuous primitive integral.  
For this we need some notation
for distributions.  The space of {\it test functions} is
$\D=C^\infty_c(\R)=\{\phi\fn\R\to\R\mid \phi\in C^\infty(\R) \text{ and }
{\rm supp}(\phi) \text{ is compact}\}$.  The {\it support} of 
function $\phi$ is the
closure of the set on which $\phi$ does not vanish and is denoted
${\rm supp}(\phi)$.  Under usual pointwise operations $\D$ is 
a linear space over field $\R$.  In $\D$ we have
a notion of convergence.  If $\{\phi_n\}\subset\D$ then $\phi_n\to 0$ as
$n\to\infty$
if there is a compact set $K\subset\R$ such that for each $n$,
${\rm supp}(\phi_n)\subset K$, and for each $m\geq 0$ we have
$\phi_n^{(m)}\to 0$ uniformly on $K$ as $n\to\infty$.  The {\it distributions}
are denoted $\Dp$ and
are the continuous linear functionals on $\D$.  For $T\in\Dp$ and
$\phi\in\D$ we write $\langle T,\phi\rangle\in\R$.  For $\phi,\psi\in\D$
and $a,b\in\R$ we have $\langle T,a\phi+b\psi\rangle=a\langle T,\phi\rangle
+ b\langle T,\psi\rangle$.  And, if $\phi_n\to 0$ in $\D$ then
$\langle T, \phi_n\rangle\to 0$ in $\R$.  Linear operations are defined
in $\Dp$ by $\langle aS+bT,\phi\rangle=a\langle S,\phi\rangle
+ b\langle T,\phi\rangle$ for $S,T\in\Dp$; $a,b\in\R$ and $\phi\in\D$.
If $f\in L^1_{\it loc}$
then $\langle T_f, \phi\rangle = \intinf f(x)\phi(x)\,dx$ defines a distribution
$T_f\in\Dp$.  The integral exists as a Lebesgue integral.
All distributions have derivatives of
all orders that are themselves distributions.
For $T\in\Dp$ and
$\phi\in\D$ the {\it distributional derivative} of $T$ is $T'$
where $\langle T',\phi\rangle = - \langle T,\phi'\rangle$.
This is also called the {\it weak derivative}.
If $p\fn\R\to\R$ is a function that is differentiable in
the pointwise sense at $x\in\R$ then we write its derivative
as $p'(x)$.
If $p$ is
a $C^\infty$ bijection such that
$p'(x)\not=0$ for any $x\in\R$  then the composition with
distribution $T$ is defined by
$\langle T\circ p,\phi\rangle =
\langle T,\frac{\phi\circ p^{-1}}{p'\circ p^{-1}}\rangle$ for all
$\phi\in\D$. 
Translations are a special case.
For $x\in\R$ define the {\it translation} $\tau_x$ on distribution
$T\in\Dp$ by $\langle \tau_xT, \phi\rangle = \langle T,\tau_{-x}\phi\rangle$
for test function $\phi\in\D$ where $\tau_x\phi(y)=\phi(y-x)$.
All of the results on distributions we use can be
found in \cite{friedlander}.

The following Banach space
will be of importance: $\balexc=\{F\fn\R\to\R\mid F\in C^0(\R),
F(-\infty)=0, F(\infty)\in\R\}$.  We use  the notation
$F(-\infty)=\lim_{x\to-\infty}F(x)$ and
$F(\infty)=\lim_{x\to\infty}F(x)$. The extended real line is
denoted $\Rbar=[-\infty,\infty]$.  The space $\balexc$ then
consists of  functions continuous on $\Rbar$ with a limit of
$0$ at $-\infty$.  We denote the functions that are continuous
on $\R$ that have real limits at $\pm\infty$ by $C^0(\Rbar)$.
Hence, $\balexc$ is properly contained in $C^0(\Rbar)$, which
is itself properly contained in the
space of uniformly
continuous functions on $\R$.
The space $\balexc$ is a Banach space under the uniform norm; $\|F\|_\infty=
\sup_{x\in\R}|F(x)|=\max_{x\in\Rbar}|F(x)|$ for $F\in\balexc$.
The {\it continuous primitive integral} is
defined by taking $\balexc$ as  the space of primitives.  The
space of integrable distributions is $\alexc=\{f\in\Dp\mid
f=F' \text{ for } F\in\balexc\}$.  If $f\in\alexc$ then
$\int_a^bf=F(b)-F(a)$ for $a,b\in\Rbar$.  The distributional
differential equation  $T'=0$ has  only constant solutions so  the
primitive  $F\in\balexc$ satisfying $F'=f$ is unique.
Integrable distributions are then
tempered and of order one.  This integral, including a discussion
of extensions to $\R^n$, is described in \cite{talviladenjoy}.
A more general integral is obtained by taking the primitives to
be regulated functions, i.e., functions with a left and right
limit at each point.  See \cite{talvilarpi}.

Examples of distributions in $\alexc$ are $T_f$ for functions $f$
that have a finite Lebesgue, Henstock--Kurzweil or wide Denjoy integral.
We identify function $f$ with the distribution $T_f$. 
Pointwise function values can be recovered from $T_f$ at points
of continuity of $f$ by evaluating
the limit $\langle T_f, \phi_n\rangle$ for a {\it delta  sequence}
converging to $x\in\R$.
This is a sequence of test functions $\{\phi_n\}\subset\D$
such that for each $n$, $\phi_n\geq 0$,
$\int_{-\infty}^\infty \phi_n=1$, and the support of $\phi_n$ tends to $\{x\}$
as $n\to\infty$. Note that if $F\in C^0(\Rbar)$ is an increasing function
with $F'(x)=0$ for almost all $x\in\R$ then the Lebesgue integral
$\int_a^b F'(x)\,dx=0$ but $F'\in\alexc$ and $\int_a^b F'=F(b)-F(a)$.
For another example of a distribution in $\alexc$, let $F\in C^0(\Rbar)$
be continuous and nowhere differentiable in the pointwise sense.
Then $F'\in\alexc$ and $\int_a^bF'=F(b)-F(a)$ for all $a,b\in\Rbar$.

The space $\alexc$ is a
Banach space under the {\it Alexiewicz norm};
$\|f\|=\sup_{I\subset\R}|\int_If|$ where  the supremum is taken
over all intervals $I\subset\R$.
An equivalent norm is $\|f\|'=\sup_{x\in\R}|\int_{-\infty}^xf|$.
The continuous primitive
integral contains the Lebesgue, Henstock--Kurzweil
and wide Denjoy integrals since their primitives 
are continuous functions.  These three spaces of functions are
not complete under the Alexiewicz norm and in fact $\alexc$
is their completion.  The lack of a Banach space has hampered
application of the Henstock--Kurzweil  integral to problems outside
of real analysis.  As we will see below, the Banach space $\alexc$
is a suitable setting  for applications of nonabsolute
integration.

We will also  need to use functions of bounded variation.  Let
$g\fn\R\to\R$.  The {\it  variation} of $g$ is 
$V\!g=\sup\sum|g(x_i)-g(y_i)|$ where the  supremum
is taken over all disjoint intervals  $\{(x_i,y_i)\}$.
The functions of  {\it bounded variation} are  denoted
$\bv=\{g\fn\R\to\R\mid V\!g<\infty\}$.  This is a Banach
space under the norm $\|g\|_\bv=|g(-\infty)|+V\!g$.
Equivalent norms are $\|g\|_\infty +Vg$ and $|g(a)|+V\!g$
for each $a\in\Rbar$.
Functions of bounded variation have a left and right limit
at each point in $\R$ and limits  at $\pm\infty$ so as above we will define 
$g(\pm\infty)=\lim_{x\to\pm\infty}g(x)$.

If $g\in L^1_{loc}$ then the {\it 
essential variation} of $g$ is ${\rm ess\,var}g=\sup\intinf g\phi'$
where the supremum is taken over  all $\phi\in\D$ with $\|\phi\|_\infty\leq 1$.
Then $\ebv=\{g\in L^1_{loc}\mid {\rm ess\,var}g<\infty\}$.
This is a Banach space under the norm $\|g\|_\ebv={\rm ess\,sup}|g|
+{\rm ess\,var}g$.
Let $0\leq\gamma\leq 1$. For $g\fn\R\to\R$ define
$g_\gamma(x)=(1-\gamma)g(x-)+\gamma g(x+)$.  For left continuity,
$\gamma=0$ and for right continuity $\gamma=1$.  The functions of
{\it normalised bounded variation} are $\nbv_\gamma=
\{g_\gamma\mid g\in\bv\}$.  If $g\in\ebv$ then
${ess\,var}g=\inf Vh$ such that $h=g$ almost everywhere. For 
each $0\leq\gamma\leq 1$
there is exactly one function
$h\in\nbv_\gamma$ such that $g=h$ almost everywhere.  In this
case ${ess\,var}g=Vh$.  Changing $g$ on a set of measure zero
does not affect its essential variation.  Each function of 
essential bounded variation has a distributional derivative that is
a signed Radon measure.  This will be denoted $\mu_g$ where
$\langle g',\phi\rangle=-\langle g,\phi'\rangle=-\intinf g\phi'=
\intinf \phi\,d\mu_g$ for all $\phi\in\D$.  

We will see that  $\ast\fn \alexc\times \bv\to C^0(\Rbar)$ and that
$\|f\ast g\|_\infty\leq \|f\|\|g\|_{\bv}$.
Similarly for $g\in\ebv$.  Convolutions  for
$f\in\alexc$ and  $g\in L^1$ will be defined using  sequences
in $\bv\cap L^1$ that  converge to $g$ in the $L^1$ norm.  It  will
be shown  that $\ast\fn \alexc\times L^1\to\alexc$ and that
$\|f\ast g\|\leq \|f\| \|g\|_1$.

Convolutions can be defined for distributions in several different
ways.
\begin{defn}\label{defndistconv}
Let $S,T\in\Dp$ and $\phi,\psi\in\D$.  Define  ${\tilde \phi}(x)=\phi(-x)$.  
(i) $\langle T\ast\psi,\phi\rangle
=\langle T,\phi\ast \tilde{\psi}\rangle$.
(ii) For each $x\in\R$, let $T\ast \psi(x)=\langle T,\tau_x\tilde{\psi}\rangle$.
(iii) $\langle S\ast T,\phi\rangle=\langle S(x),\langle T(y),\phi(x+y)
\rangle\rangle$.
\end{defn}
In (i), $\ast\fn \Dp\times\D\to\Dp$.  This definition also applies to
other
spaces of  test functions and their  duals, such as the Schwartz
space of rapidly decreasing functions or the compactly supported
distributions.
In (ii), $\ast\fn\Dp\times\D\to C^\infty$.
In \cite{folland} it is shown that definitions (i) and (ii) are equivalent.
In (iii), $\ast\fn \Dp\times\Dp\to\Dp$.  However, this definition requires
restrictions on the supports of $S$ and $T$.  It suffices that one of
these distributions have compact support.  Other conditions on the supports
can be imposed.  See \cite{friedlander} and \cite{zemanian}.
This definition is an instance of the tensor product,
$\langle S\otimes T,\Phi\rangle=\langle S(x),\langle T(y),\Phi(x,y)\rangle
\rangle$ where now $\Phi\in\D(\R^2)$.

Under (i), $T\ast \psi$ 
is in $C^\infty$.  It
satisfies $(T\ast\psi)\ast\phi=T\ast(\psi\ast\phi)$, 
$\tau_x(T\ast \psi)=(\tau_xT)\ast\psi=T\ast(\tau_x\psi)$,
and
$(T\ast\psi)^{(n)}=T\ast\psi^{(n)}
=T^{(n)}\ast\psi$.
Under (iii), with appropriate support restrictions,
$S\ast T$ is in $\Dp$.  It is commutative and 
associative, commutes with translations, and satisfies
$(S\ast T)^{(n)}=S^{(n)}\ast T=S\ast T^{(n)}$.  It is weakly
continuous in $\Dp$, i.e., if $T_n\to T$ in $\Dp$ then
$T_n\ast\psi\to T\ast \psi$ in $\Dp$.  See \cite{folland},
\cite{friedlander}, \cite{reedsimon} and \cite{zemanian}
for additional properties of convolutions of distributions.

Although elements of $\alexc$ are  distributions, we show in
this paper that their behaviour
as convolutions is more like that of integrable functions.

An appendix contains the proof of a type of Fubini theorem.

\section{Convolution in  $\alexc\times\bv$}
In this section we prove basic results for the convolution
when  $f\in\alexc$ and $g\in\bv$.  Under these
conditions $f\ast g$ is commutative, continuous on $\Rbar$ and
commutes with translations.
It can be estimated in the uniform norm in terms of the Alexiewicz
and $\bv$ norms.  There is also an associative property.
We first need
the result that $\bv$ forms the space of multipliers for
$\alexc$, i.e., if  $f\in\alexc$ then $fg\in\alexc$ for all
$g\in\bv$.  The  integral $\int_{I}fg$ is  defined using the
integration by parts formula in the Appendix.  The H\"older
inequality \eqref{holder} shows that $\bv$ is the dual space of $\alexc$.

We define the convolution of $f\in\alexc$ and $g\in\bv$ as
$f\ast g(x)=\intinf (f\circ r_x)g$ where $r_x(t)=x-t$.  We
write this as $f\ast g(x)=\intinf f(x-y)g(y)\,dy$.
\begin{theorem}\label{ginbv}
Let $f\in\alexc$ 
and let $g\in\bv$. Then (a) $f\ast g$  exists on $\R$
(b) $f\ast g=g\ast f$  (c) $\|f\ast g\|_\infty \leq
|\int_{-\infty}^\infty f|\inf_{\R}|g| + \|f\|Vg\leq \|f\| \|g\|_\bv$
(d) $f\ast g\in C^0(\Rbar)$,  $\lim_{x\to\pm\infty}
f\ast g(x)=g(\pm\infty)\intinf f$.
(e) If $h\in L^1$ then $f\ast (g\ast h)=(f\ast g)\ast h\in C^0(\Rbar)$.
(f) Let $x,z\in\R$.  Then 
$\tau_z(f\ast g)(x)=(\tau_zf)\ast g(x)=(f\ast \tau_zg)(x)$.
(g) For each $f\in\alexc$ define $\Phi_f\fn \bv\to C^0(\Rbar)$ by
$\Phi_f[g]=f\ast g$.  Then $\Phi_f$ is a bounded linear operator
and $\|\Phi_f\|\leq \|f\|$. There exists a nonzero distribution $f\in\alexc$
such that $\|\Phi_f\|=\|f\|$.   For each $g\in\bv$ define 
$\Psi_g\fn \alexc\to C^0(\Rbar)$ by
$\Psi_g[f]=f\ast g$.  Then $\Psi_g$ is a bounded linear operator
and $\|\Psi_g\|\leq\|g\|_\bv$. 
There exists a nonzero function
$g\in\bv$ such that
$\|\Psi_g\|=\|g\|_{\bv}$. 
(h) ${\rm supp}(f\ast g)\subset {\rm supp}(f)+{\rm supp}(g)$.
\end{theorem}
\bigskip
\noindent
{\bf Proof:}  (a) Existence is given via the integration by parts
formula \eqref{parts} in the Appendix.  (b) See \cite[Theorem~11]{talviladenjoy}
for a change of variables theorem that can be  used with $y\mapsto x-y$.
(c) This inequality follows from the H\"older inequality \eqref{holder}.
(d) Let $x,t\in\R$.  From (c) we have
\begin{eqnarray*}
|f\ast g(t)-f\ast g(x)| & \leq & \|f(t-\cdot) - f(x-\cdot)\| \|g\|_\bv\\
  & = & \|f(t-x-\cdot) - f(\cdot)\| \|g\|_\bv\\
 & \to & 0 \text{ as } t\to x.
\end{eqnarray*}
The last line follows from continuity in the  Alexiewicz norm 
\cite[Theorem~22]{talviladenjoy}.  Hence, $f\ast g$  is  uniformly continuous
on $\R$.  And, $\lim_{x\to\infty}\intinf f(y)g(x-y)\,dy
=\intinf f(y)\lim_{x\to\infty}g(x-y)\,dy=g(\infty)\intinf  f$.
The limit $x\to\infty$  can be taken under the  integral sign
since $g(x-y)$ is of  uniform bounded variation, i.e., 
$V_{y\in\R}g(x-y)=V\!g$.  Theorem~22 in \cite{talviladenjoy} then applies.
Similarly as $x\to-\infty$.
(e) First show $g\ast h\in\bv$.  Let $\{(s_i,t_i)\}$ be disjoint intervals
in $\R$.  Then
\begin{eqnarray*}
\sum|g\ast h(s_i)-g\ast h(t_i)| & \leq &  
\sum\intinf|g(s_i-y)-g(t_i-y)| |h(y)|\,dy\\
 & = & \intinf\sum|g(s_i-y)-g(t_i-y)| |h(y)|\,dy.
\end{eqnarray*}
Hence, $V(g\ast h)\leq V\!g\, \|h\|_1$.  The interchange of sum and
integral follows from the Fubini--Tonelli theorem.  Now (d) shows
$f\ast (g\ast h)\in C^0(\Rbar)$.  Write
\begin{eqnarray*}
f\ast (g\ast h)(x) & = & \intinf f(y)\intinf g(x-y-z)h(z)\,dz\,dy\\
 & = & \intinf h(z)\intinf f(y)g(x-y-z)\,dy\,dz\\
 & = & (f\ast g)\ast h(x).
\end{eqnarray*}
We can interchange orders of integration using Proposition~\ref{fubini}.
For (ii) in Proposition~\ref{fubini}, the function 
$z\mapsto V_{y\in\R}g(x-y-z)h(z)=V\!g\, h(z)$ is in $L^1$
for each fixed $x\in\R$.  Since $g$ is of bounded variation it is bounded
so $|g(x-y-z)h(z)|\leq \|g\|_\infty |h(z)|$ and condition (iii) is
satisfied. (f) This follows from a linear change of variables as in (a).
(g) From  (c) we have $\|\Phi_f\|=\sup_{\|g\|_\bv=1}\|f\ast g\|_\infty
\leq  \sup_{\|g\|_\bv=1}\|f\|\|g\|_\bv=\|f\|$.
Let $f>0$ be in $L^1$.  If $g=1$ then $\|g\|_\bv=1$ and
$f\ast g(x)=\intinf f$ so $\|\Phi_f\|=\|f\|=\|f\|_1$.
To prove $\|\Psi_g\|\leq \|g\|_\bv$, note
that $\|\Psi_g\|=\sup_{\|f\|=1}\|f\ast g\|_\infty
\leq  \sup_{\|f\|=1}\|f\|\|g\|_\bv=\|g\|_\bv$.
Let $g=\chi_{(0,\infty)}$.
Then
$\|\Psi_g\|=\sup_{\|f\|=1} \|f\ast g\|_\infty =\sup_{\|f\|=1}\sup_{x\in\R}
|\int_{-\infty}^xf|= 1=
\|g\|_\bv$.
(h) Suppose $x\notin {\rm supp}(f)+{\rm supp}(g)$.
Note that we can write $f\ast g(x)=\intinf g(x-y)\,dF(y)$, in terms
of a Henstock--Stieltjes integral.  See \cite{talviladenjoy} for details.
This integral is approximated by Riemann sums $\sum_{n=1}^N g(x-z_n)
[F(t_n)-F(t_{n-1})]$ where $z_n\in[t_{n-1},t_n]$,
$-\infty=t_0<t_1<\ldots<t_N=\infty$ and there is a gauge function
$\gamma$ mapping $\Rbar$ to the open intervals in $\Rbar$ such that
$[t_{n-1},t_n]\subset\gamma(z_n)$.  If $z_n\notin \supp(f)$ then
since $\R\setminus\supp(f)$ is open there is an open interval
$z_n\subset I\subset \R\setminus\supp(f)$.  We can take $\gamma$
such that $[t_{n-1},t_n]\subset I$ for all $1\leq n\leq N$.  And,
$F$ is constant on each interval in $\R\setminus\supp(f)$.  Therefore,
$g(x-z_n)[F(t_n)-F(t_{n-1})]=0$ and only tags $z_n\in\supp(f)$ can
contribute to the Riemann sum.  But for all $z_n\in\supp(f)$ we have
$x-z_n\notin\supp(g)$ so $g(x-z_n)[F(t_n)-F(t_{n-1})]=0$.  It follows
that $f\ast g(x)=0$.
\qed

Similar results are proven for $f\in L^p$ in \cite[\S~8.2]{folland}.

If we use the equivalent norm $\|f\|'=\sup_{x\in\R}|\int_{-\infty}^xf|$
then $\|\Phi_f\|=\|f\|'$.   For, integration by parts gives
$\|\Phi_f\|\leq \|f\|'$.
Now given $f\in\alexc$ let $g=\chi_{(0,\infty)}$.  Then
$\|g\|_\bv=1$.  And, $f\ast g(x)=
\int_{-\infty}^xf$.  Hence, $\|f\ast g\|_\infty=\|f\|'$ and
$\|\Phi_f\|=\|f\|'$. 
We can have strict inequality in $\|\Psi_g\|\leq \|g\|_\bv$.  For
example, let $g=\chi_{\{0\}}$.  Then $\|g\|_\bv=2$ but integration
by parts shows $f\ast g=0$ for each $f\in\alexc$.

\begin{remark}
{\rm
If $f\in\alexc$ and $g\in\ebv$ we can use Definition~\ref{defnebv} to
define $f\ast g(x)=f\ast g_\gamma(x)$ where $g_\gamma=g$ almost
everywhere and $g_\gamma\in\nbv_\gamma$.  All of the results in
Theorem~\ref{ginbv} and the rest of the paper have analogues.
Note that $f\ast g(x)=F(\infty)g_\gamma(-\infty)+F\ast \mu_g$.
}
\end{remark}

\begin{prop}
The three definitions of convolution for distributions in 
Definition~\ref{defndistconv} are compatible with $f\ast g$
for $f\in\alexc$ and $g\in\bv$.
\end{prop}
\bigskip
\noindent
{\bf Proof:} Let $f\in\alexc$, $g\in\bv$ and $\phi,\psi\in\D$.
Definition~\ref{defndistconv}(i) gives
\begin{eqnarray*}
\langle f,{\tilde \psi}\ast\phi\rangle & = & \intinf f(x)\intinf \psi(y-x)\phi(y)
\,dy\,dx\\
 & = & \intinf \intinf f(x)\psi(y-x)\phi(y)\,dx\,dy\\
 & = & \langle f\ast\psi,\phi\rangle.
\end{eqnarray*}
Since $\psi\in\bv$ and $\phi\in L^1$, Proposition~\ref{fubini}
justifies the interchange of integrals.
Definition~\ref{defndistconv}(ii) gives 
$$
\langle f,\tau_x{\tilde \psi}\rangle  =  \intinf f(y)\psi(x-y)\,dy=
f\ast\psi(x).
$$
Definition~\ref{defndistconv}(iii) gives 
\begin{eqnarray*}
\langle f(y),\langle g(x),\phi(x+y)\rangle\rangle & = &
\intinf f(y)\intinf g(x)\phi(x+y)\,dx\,dy\\
 & = & \intinf f(y)\intinf g(x-y)\phi(x)\,dx\,dy\\
 & = & \intinf \phi(x)\intinf f(y) g(x-y)\,dx\,dx\\
 & = & \langle f\ast g,\phi\rangle.
\end{eqnarray*}
The interchange of integrals is accomplished using
Proposition~\ref{fubini} since $g\in\bv$ and $\phi\in L^1$. \qed

The locally integrable distributions are defined as
$\alexc(loc)=\{f\in\Dp\mid f=F' \text{ for some } F\in C^0(\R)\}$.
Let $f\in\alexc(loc)$ and let $g\in\bv$ with support in the compact
interval $[a,b]$.  By the Hake theorem \cite[Theorem~25]{talviladenjoy},
$f\ast g(x)$ exists if and only if  the limits of $\int_{\alpha}^\beta
f(x-y)g(y)\,dy$ exist as $\alpha\to-\infty$ and $\beta\to\infty$.
This gives
$$
f\ast g(x)  =  \int_a^b f(x-y)g(y)\,dy =\int_{x-b}^{x-a} f(y)g(x-y)\,dy.
$$
There are analogues of the results in Theorem~\ref{ginbv}.  For
example, $|f\ast g(x)|\leq |\int_{x-b}^{x-a} f|\inf_{[a,b]}|g|
+ \|f\chi_{[x-b,x-a]}\|V_{[a,b]}g$.  There are also versions where
the supports are taken to be semi-infinite intervals.

We can also define the distributions with bounded primitive as
$\alexc(bd)=\{f\in\Dp\mid f=F' \text{ for some bounded } F\in C^0(\R)
\text{ with } F(0)=0\}$.
Let $f\in\alexc(bd)$ and let $F$ be its unique primitive.
If $g\in\bv$ such that $g(\pm\infty)=0$ then
\begin{eqnarray*}
f\ast g(x) & = & 
\llim_{\stackrel{\alpha\to-\infty}{{\scriptscriptstyle\beta\to\infty}}}
\int_\alpha^\beta f(x-y)g(y)\,dy\\
 & = & \llim_{\stackrel{\alpha\to-\infty}{{\scriptscriptstyle\beta\to\infty}}}\left[
F(x-\alpha)g(\alpha)-F(x-\beta)g(\beta)+\int_\alpha^\beta F(x-y)\,dg(y)\right
]\\
 & = & \intinf F(x-y)\,dg(y) = \intinf F(y)\,dg(x-y).
\end{eqnarray*}
It follows that $\|f\ast g\|_\infty\leq \|F\|_\infty Vg$.

It is possible to formulate other existence criteria.  For example,
if $f(x)=\log|x|\sin(x)$ and $g(x)=|x|^{-\alpha}$ for some
$0<\alpha<1$ then $f$ and $g$ are not in $\alexc$, $\bv$ or
$L^p$ for any $1\leq p\leq\infty$ but $f\ast g$ exists on
$\R$ because $f,g\in L^1_{loc}$ and if $F(x)=\int_0^xf$ then
$\lim_{|x|\to\infty}F(x)g(x)=0$.

The following example shows that $f\ast g$ need not be of bounded variation
and hence not absolutely continuous.
Let $g=\chi_{(0,\infty)}$.  For
$f\in\alexc$ we have
$
f\ast g(x)=\int_{-\infty}^x f=F(x)$
where $F\in\balexc$ is the primitive of $f$. But $F$ need not be  
of bounded variation or even of  local bounded variation.  For example, 
let $f(x)=\sin(x^{-2})-2x^{-2}\cos(x^{-2})$ and  let  
$F$ be its primitive in $\balexc$.
Finally, although $f\ast g$ is continuous it need not be  integrable over $\R$.
For example, let $g=1$ then $f\ast g(x)=\intinf f$ and $\intinf f\ast g$
only exists  if $\intinf  f=0$.

\section{Convolution in $\alexc\times L^1$}
We now extend the convolution $f\ast g$ to $f\in\alexc$ and
$g\in L^1$.  Since there are functions in $L^1$ that are not of bounded
variation, there are distributions $f\in \alexc$ and functions
$g\in L^1$ such that the integral $\intinf f(x-y)g(y)\,dy$ does not
exist. The convolution is then defined as the limit in $\|\cdot\|$
of a sequence
$f\ast g_n$ for $g_n\in\bv\cap L^1$ such that $g_n\to g$ 
in the $L^1$ norm.  This is possible since $\bv\cap L^1$ is dense
in $L^1$.  We also give an equivalent definition using
the fact that $L^1$ is dense in $\alexc$.  Take a sequence $\{f_n\}\subset
L^1$ such that $\|f_n-f\|\to 0$.  Then $f\ast g$ is the 
limit in $\|\cdot\|$ of $f_n\ast g$.  In this more general setting of
convolution defined in $\alexc\times L^1$ we now have
an Alexiewicz norm estimate for $f\ast g$ in terms of estimates
of $f$ in the Alexiewicz norm and $g$ in the $L^1$ norm.  There
is associativity with $L^1$ functions and commutativity with translations.
\begin{defn}\label{L1defn}
Let $f\in\alexc$  and let $g\in L^1$.  Let $\{g_n\}\subset \bv\cap L^1$
such that $\|g_n -g\|_1\to 0$.  Define $f\ast g$ 
as the unique element
in $\alexc$ such that $\|f\ast g_n -f\ast g\|\to 0$.
\end{defn}
To see that the definition makes sense, first note that $\bv\cap L^1$ is
dense in $L^1$ since step functions are dense in $L^1$.  Hence, the
required sequence $\{g_n\}$ exists.  Let $[\alpha,\beta]\subset\R$ be a  compact interval.
Let $F\in\balexc$ be the primitive  of  $f$.
Then
\begin{eqnarray}
\int_\alpha^\beta f\ast g_n(x)\,dx & = & \int_\alpha^\beta\intinf f(y)g_n(x-y)
\,dy\,dx\notag\\
 & = & \intinf f(y)\int_\alpha^\beta g_n(x-y)\,dx\,dy\label{theorem1a}\\
 & = &  \intinf f(y)\int_{\alpha-y}^{\beta-y} g_n(x)\,dx\,dy\notag\\
 & = & -\intinf F(y)\, d\left[\int_{\alpha-y}^{\beta-y}g_n\right]\label{theorem1b}\\
 & = & \intinf F(y)\left[g_n(\beta-y)-g_n(\alpha-y)\right]dy\label{theorem1c}\\
 & = & \intinf \left(\int_{\alpha-y}^{\beta-y}f\right)g_n(y)\,dy.\notag
\end{eqnarray}
The interchange of orders of integration in \eqref{theorem1a}  
is accomplished with
Proposition~\ref{fubini} using $g(x,y)=g_n(x-y)\chi_{[\alpha,\beta]}(x)$.
Integration by parts  gives \eqref{theorem1b} since
$\lim_{y\to\infty}\int_{\alpha-y}^{\beta-y}g_n=0$.
As $F$ is continuous  and the function $y\mapsto \int_{\alpha-y}^{\beta-y}g_n$
is absolutely continuous we  get \eqref{theorem1c}.  Taking the
supremum over $\alpha,\beta\in\R$ gives 
\begin{equation}
\|f\ast g_n\|\leq \|f\| \|g_n\|_1.\label{bvL1inequality}
\end{equation}
We now have
$$
\|f\ast g_m-f\ast g_n\|  =  \|f\ast(g_m-g_n)\|
 \leq  \|f\|\|g_m-g_n\|_1
$$
and $\{f\ast g_n\}$ is a Cauchy sequence in $\alexc$. Since $\alexc$ 
is complete this sequence has a limit in $\alexc$ which we denote
$f\ast g$.  The definition does not depend on the choice of sequence
$\{g_n\}$, for if $\{h_n\}\subset\bv\cap L^1$ such that $\|h_n-g\|_1\to 0$
then $\|f\ast g_n-f\ast h_n\|\leq \|f\|(\|g_n-g\|_1+\|h_n-g\|_1)\to 0$
as $n\to\infty$.  The above calculation also shows that if $g\in\bv\cap L^1$
then the integral definition $f\ast g(x)=\intinf f(x-y)g(y)\,dy$ and the
limit definition agree.

\begin{defn}\label{L1defn2}
Let $f\in\alexc$  and let $g\in L^1$.  Let $\{f_n\}\subset L^1$
such that $\|f_n -f\|\to 0$.  Define $f\ast g$ 
as the unique element
in $\alexc$ such that $\|f_n\ast g -f\ast g\|\to 0$.
\end{defn}
To show this definition makes sense, first show $L^1$ is
dense in $\alexc$.

\begin{prop}\label{acr}
$L^1$  is dense in $\alexc$.
\end{prop}
\bigskip
\noindent
{\bf Proof:}
Let $AC(\Rbar)$ be the functions that are absolutely continuous
on each compact interval and which are of bounded variation on
the real line.
Then $f\in L^1$ if and only if there exists $F\in\acr$ such that
$F'(x)=f(x)$ for almost all $x\in\R$.  
Let $f\in\alexc$ be given.
Let $F\in\balexc$ be its primitive. For $\epsilon>0$, take
$M>0$ such that $|F(x)|<\epsilon$ for $x<-M$ and $|F(x)-F(\infty)|
<\epsilon$  for $x>M$.  Due to the Weierstrass approximation theorem
there is a continuous function $P\fn\R\to\R$ such that $P(x)=F(-M)$ for $x\leq -M$, 
$P(x)=F(M)$ for $x\geq M$,
$|P(x)-F(x)|<\epsilon$ for $|x|\leq M$ and $P$ is a polynomial on $[-M,M]$.
Hence, $P\in\acr$ and $\|P'-f\|<3\epsilon$.
\qed

In Definition~\ref{L1defn2}, the required sequence
$\{f_n\}\subset L^1$ exists.  
Let $[\alpha,\beta]\subset\R$ be a  compact interval.
Then, by the usual Fubini--Tonelli theorem in $L^1$,
\begin{eqnarray*}
\int_\alpha^\beta f_n\ast g(x)\,dx & = & \int_\alpha^\beta\intinf f_n(x-y)g(y)
\,dy\,dx\notag\\
 & = & \intinf g(y)\int_\alpha^\beta f_n(x-y)\,dx\,dy.\\
\end{eqnarray*}
Take the
supremum over $\alpha,\beta\in\R$ and use the
$L^1-L^\infty$ H\"older inequality to  get
\begin{equation}
\|f_n\ast g\|\leq \|f_n\| \|g\|_1.\label{L1defn2inequality}
\end{equation}
It now follows that
$\{f_n\ast g\}$ is a Cauchy sequence.  It then converges to
an element of $\alexc$.  Inequality \eqref{L1defn2inequality} 
also shows this limit
is independent of the choice of $\{f_n\}$.  To see that Definition~\ref{L1defn}
and  Definition~\ref{L1defn2} agree, take $\{f_n\}\subset L^1$ with
$\|f_n-f\|\to 0$  and
$\{g_n\}\subset \bv\cap L^1$ with $\|g_n-g\|_1\to 0$.  Then
\begin{eqnarray*}
\|f_n\ast g-f\ast g_n\| & = & \|(f_n-f)\ast g -f\ast(g_n-g)\|\\
 & \leq &  \|(f_n-f)\ast g\| + \|f\ast(g_n-g)\|\\
 & \leq &  \|f_n-f\|\|g\|_1 + \|f\|\|g_n-g\|_1.
\end{eqnarray*}
Letting $n\to\infty$ shows the limits of $f_n\ast g$ in Definition~\ref{L1defn2}
and $f\ast g_n$ in Definition~\ref{L1defn} are the same.

\begin{theorem}\label{L1convolution}
Let $f\in\alexc$ and $g\in L^1$.  Define $f\ast g$ as in 
Definition~\ref{L1defn}.  Then (a) $\|f\ast g\|\leq \|f\| \|g\|_1$.
(b) Let $h\in L^1$.  Then $(f\ast g)\ast h=f\ast (g\ast h)\in\alexc$.
(c) For each $z\in\R$, $\tau_z(f\ast g)=(\tau_zf)\ast g=(f\ast \tau_zg)$.
(d)
For each $f\in\alexc$ define $\Phi_f\fn L^1\to \alexc$ by
$\Phi_f[g]=f\ast g$.  Then $\Phi_f$ is a bounded linear operator
and $\|\Phi_f\|\leq\|f\|$.   There exists a nonzero distribution
$f\in\alexc$ such that
$\|\Phi_f\|=\|f\|$.  For each $g\in L^1$ define 
$\Psi_g\fn \alexc\to \alexc$ by
$\Psi_g[f]=f\ast g$.  Then $\Psi_g$ is a bounded linear operator
and $\|\Psi_g\|\leq\|g\|_1$.
There exists a nonzero function
$g\in L^1$ such that
$\|\Psi_g\|=\|g\|_{\bv}$.
(e) Define $g_t(x)=g(x/t)/t$ for $t>0$.
Let $a=\intinf g_t(x)\,dx=\intinf g$.  Then $\|f\ast g_t-af\|\to 0$
as $t\to 0$.
(f)  ${\rm supp}(f\ast g)\subset {\rm supp}(f)+{\rm supp}(g)$.
\end{theorem}
\bigskip
\noindent
{\bf Proof:} Let $\{g_n\}$ be as in Definition~\ref{L1defn}.
(a)
Since $\|f\ast g_n\|\to\|f\ast g\|$, equation \eqref{bvL1inequality}
shows 
$\|f\ast g\|\leq \|f\| \|g\|_1$.
(b) 
Let $\{h_n\}
\subset\bv\cap L^1$ such that $\|h_n-h\|_1\to 0$.  Then
$(f\ast g)\ast h:=\xi\in\alexc$ such that $\|(f\ast g)\ast h_n-\xi\|
\to 0$.  Since $g\ast h\in L^1$ there is $\{p_n\}
\subset\bv\cap L^1$ such that $\|p_n-g\ast h\|_1\to 0$.
Then $f\ast( g\ast h):=\eta\in\alexc$ such that $\|f\ast p_n-\eta\|
\to 0$.  Now,
\begin{eqnarray*}
\|\xi-\eta\| & \leq & \|(f\ast g)\ast h_n-\xi\| + \|f\ast p_n-\eta\|\\
 & & \quad + \|(f\ast g)\ast h_n-(f\ast g_n)\ast h_n\| 
+ \|(f\ast g_n)\ast h_n-f\ast p_n\|.
\end{eqnarray*}
Using \eqref{bvL1inequality},
\begin{eqnarray*}
\|(f\ast g)\ast h_n-(f\ast g_n)\ast h_n\| & = &
\|[f\ast(g-g_n)]\ast h_n\|\\
 & \leq & \|f\| \|g_n-g\|_1 \|h_n\|_1\\
 & \to & 0 \text{ as } n\to \infty.
\end{eqnarray*}
Finally, use Theorem~\ref{ginbv}(e) and \eqref{bvL1inequality} to write
\begin{eqnarray*}
\|(f\ast g_n)\ast h_n-f\ast p_n\| & = & \|f\ast(g_n\ast h_n-p_n)\|\\
 & \leq & \|f\| (\|g_n-g\|_1 \|h_n\|_1 + \|g\|_1 \|h_n-h\|_1\\
 & & \qquad + 
\|p_n-g\ast h\|_1)\\
 & \to & 0 \text{ as } n\to \infty.
\end{eqnarray*}
(c) The Alexiewicz norm is invariant under translation
\cite[Theorem~28]{talviladenjoy} so $\tau_z(f\ast g)\in\alexc$. 
Use Theorem~\ref{ginbv}(f) to write
$\|\tau_z(f\ast g)- \tau_z(f\ast g_n)\|=\|f\ast g- f\ast g_n\|=
\|\tau_z(f\ast g)- (\tau_zf)\ast g_n)\|=
\|\tau_z(f\ast g)- f\ast (\tau_z g_n)\|$.
Translation invariance of the $L^1$ norm
completes the proof.
(d)  From  (a) we have $\|\Phi_f\|=\sup_{\|g\|_1=1}\|f\ast g\|
\leq  \sup_{\|g\|_1=1}\|f\|\|g\|_1=\|f\|$.
We get equality by considering $f$ and $g$ to be positive functions
in $L^1$.
To prove $\|\Psi_g\|\leq \|g\|_1$, note
that $\|\Psi_g\|=\sup_{\|f\|=1}\|f\ast g\|
\leq  \sup_{\|f\|=1}\|f\|\|g\|_1=\|g\|_1$.
We get equality by considering $f$ and $g$ to be positive functions
in $L^1$.
(e) First consider $g\in\bv\cap L^1$.
We have 
$$
f\ast g_t(x)=\intinf f(x-y)\,g\!\left(\frac{y}{t}\right)\frac{dy}{t}
=\intinf f(x-ty)g(y)\,dy.
$$
For $-\infty<\alpha<\beta<\infty$,
\begin{eqnarray}
\left|\int_\alpha^\beta\left[f\ast g_t(x)-a\,f(x)\right]dx\right|
 & = & \left|\int_\alpha^\beta\intinf\left[f(x-ty)-f(x)\right]
g(y)\,dy\,dx\right|\notag\\
 & = &  \left|\intinf\int_\alpha^\beta\left[f(x-ty)-f(x)\right]
g(y)\,dx\,dy\right|\label{gt}\\
 & \leq & \intinf\|\tau_{ty}f-f\||g(y)|\,dy\label{gt2}\\
 & \leq & 2\|f\|\|g\|_1.\notag
\end{eqnarray}
By dominated convergence we can take the limit $t\to 0$ inside
the integral \eqref{gt2}.  Continuity of $f$ in the Alexiewicz norm 
then shows $\|f\ast g_t-af\|\to 0$
as $t\to 0$.  

Now take a sequence $\{g^{(n)}\}\subset\bv\cap L^1$ such that
$\|g^{(n)}-g\|_1\to 0$.  Define $g_t^{(n)}(x)=g^{(n)}(x/t)/t$
and $a^{(n)}=\intinf g^{(n)}(x)\,dx$.  We have
\begin{equation}
\|f\ast g_t-af\|\leq \|f\ast g_t^{(n)} -a^{(n)}f\| + \|f\ast g_t^{(n)} -
f\ast g_t\|+\|a^{(n)}f-af\|.\label{gt3}
\end{equation}
By the inequality in (a), $\|f\ast g_t^{(n)} -
f\ast g_t\|  \leq  \|f\|\| g_t^{(n)} - g_t\|_1$.  Whereas,
$$
\| g_t^{(n)} - g_t\|_1 = \intinf\left|g^{(n)}\!\left(\frac{x}{t}
\right)-g\!\left(\frac{x}{t}\right)\right|\frac{dx}{t}\\
=\|g^{(n)}-g\|_1\to 0 \text{ as } n\to\infty.
$$
And, 
$
\|a^{(n)}f-af\|=|a^{(n)}-a|\|f\|=\|g^{(n)}-g\|_1\|f\|
$.
Given $\epsilon>0$ fix $n$ large enough so that
$\|f\ast g_t^{(n)} -
f\ast g_t\|+\|a^{(n)}f-af\| <\epsilon$.
Now let $t\to 0$ in \eqref{gt3}.

The interchange of order of integration in \eqref{gt}
is  justified as follows.  A change of variables and
Proposition~\ref{fubini} give
\begin{eqnarray*}
\int_\alpha^\beta\intinf f(x-ty)
g(y)\,dy\,dx
 & = & 
\int_\alpha^\beta\intinf f(y)g\!\left(\frac{x-y}{t}\right)\frac{dy}{t}\,dx\\
 & = &
\intinf\int_\alpha^\beta f(y)\,g\!\left(\frac{x-y}{t}\right)dx\,\frac{dy}{t}.\\
\end{eqnarray*}
And,
\begin{eqnarray*}
\intinf\int_\alpha^\beta f(x-ty)
g(y)\,dx\,dy
 & = & 
\intinf\intinf f(x)\,g\!\left(\frac{y}{t}\right)\chi_{(\alpha-y,\beta-y)}(x)
dx\,\frac{dy}{t}\\
 & = & 
\intinf\intinf f(x)\,g\!\left(\frac{y}{t}\right)\chi_{(\alpha-y,\beta-y)}(x)
\frac{dy}{t}\,dx\\
 & = & 
\intinf\int_\alpha^\beta f(x)\,g\!\left(\frac{y-x}{t}\right)\frac{dy}{t}\,dx.
\end{eqnarray*}
Note that $\int_\alpha^\beta\intinf f(x)g(y)\,dy\,dx=
\intinf\int_\alpha^\beta f(x)g(y)\,dx\,dy$
by Corollary~\ref{corfubini}. (f) This follows from the equivalence
of Definition~\ref{defndistconv} and Definition~\ref{L1defn}, 
proved in Proposition~\ref{distL1defn}.
See Theorem~5.4-2 and Theorem~5.3-1 in \cite{zemanian}.\qed

Young's inequality states that
$\|f\ast g\|_p\leq \|f\|_p\|g\|_1$
when $f\in L^p$ for some $1\leq p\leq \infty$ and $g\in L^1$.
Part (a) of Theorem~\ref{L1convolution} extends this to $f\in \alexc$.
See \cite{folland} for other
results when $f\in L^p$.

The fact that convolution is linear in both arguments, together with
(b), shows that $\alexc$ is an $L^1$-module over the $L^1$ convolution
algebra.  See \cite{dales} for the definition.  It does not appear
that $\alexc$ is a Banach algebra under convolution.

We now show that Definition~\ref{defndistconv}(iii) and the above
definitions agree.
\begin{prop}\label{distL1defn}
Let $f\in\alexc$, $g\in L^1$ and $\phi\in\D$.  Define
$F(y)=\int_{-\infty}^y f$ and $G(x)=\int_{-\infty}^x g$.  
Definition~\ref{defndistconv} and Definition~\ref{L1defn}
both give
\begin{eqnarray}
\langle f\ast g,\phi\rangle & = & \intinf f(y)\intinf g(x)\phi(x+y)\,dx\,dy
\label{distint1}\\
  & = &  \intinf\intinf F(y)G(x)\phi''(x+y)\,dx\,dy.\label{distint2}
\end{eqnarray}
\end{prop}
\bigskip
\noindent
{\bf Proof:}  Let $\Phi(y)=\intinf g(x)\phi(x+y)\,dx$.  Then $\Phi\in
C^\infty(\R)$ and $\Phi'(y)=\intinf g(x)\phi'(x+y)\,dx$.
And, $\intinf|\Phi'(y)|\,dy\leq\intinf|g(x)|\intinf|\phi'(x+y)|\,dy\,dx
\leq \|g\|_1\|\phi'\|_1$ so $\Phi\in \acr$.  Dominated convergence
then shows $\lim_{|y|\to\infty}\Phi(y)=0$.  Integration by parts now
gives \eqref{distint1} and \eqref{distint2}.

Let $\{g_n\}\subset \bv\cap L^1$
such that $\|g_n -g\|_1\to 0$.  Since convergence in $\|\cdot\|$
implies convergence in $\Dp$, we have
\begin{eqnarray*}
\langle f\ast g,\phi\rangle & = & \lim_{n\to\infty}
\langle f\ast g_n,\phi\rangle\\
 & = & \lim_{n\to\infty}\intinf\intinf f(y)g_n(x-y)\phi(x)\,dy\,dx\\
 & = & \lim_{n\to\infty}\intinf f(y)\intinf g_n(x-y)\phi(x)\,dx\,dy.\\
\end{eqnarray*}
Proposition~\ref{fubini} allows interchange of the iterated integrals.
Define $\Phi_n(y)=\intinf g_n(x)\phi(x+y)\,dx$.
As above, $V\Phi_n\leq \|g_n\|_1\|\phi'\|_1\leq (\|g\|_1+1)\|\phi'\|_1$
for large enough $n$.  Hence, $\Phi_n$ is of uniform bounded variation.
Theorem~22 in \cite{talviladenjoy} then gives
$\langle f\ast g,\phi\rangle =\intinf f(y)\lim_{n\to\infty}\Phi_n(y)\,dy
=\intinf f(y)\intinf g(x-y)\phi(x)\,dx\,dy$.  The last step follows
since $\|g_n-g\|_1\to 0$.\qed

If $g\in L^1\setminus\bv$ then $f\ast g$ need not be continuous
or bounded.
For example, take $1/2\leq\alpha<1$ and let 
$f(x)=g(x)=x^{-\alpha}\chi_{(0,1)}(x)$.
Then
$f\in L^1\subset\alexc$ and $g\in L^1\setminus\bv$.  We have
$f\ast g(x)=0$ for $x\leq 0$.  For $0<x\leq 1$ we have
$f\ast g(x)=\int_{0}^xy^{-\alpha}(x-y)^{-\alpha}\,dy=x^{1-2\alpha}
\int_0^1 y^{-\alpha}(1-y)^{-\alpha}
\,dy=x^{1-2\alpha}\Gamma^2(1-\alpha)/\Gamma(2-2\alpha)$.  
Hence, $f\ast g$ is not continuous at $0$.  If $1/2<\alpha<1$
then $f\ast g$ is unbounded at $0$.

As another example, consider $f(x)=\sin(\pi x)/\log|x|$ and 
$g(x)=\chi_{(0,1)}(x)$.  Then $f\in\alexc$ and for each
$1\leq p\leq\infty$ we have $g\in\bv\cap L^p$.  And,
\begin{eqnarray*}
f\ast g(x) & = & \int_{x-1}^x\frac{\sin(\pi y)}{\log(y)}dy\quad\text{ for }
x\geq 2\\
 & = & \frac{\cos(\pi(x-1))}{\pi\log(x-1)}-\frac{\cos(\pi x)}{\pi\log(x)}
- \frac{1}{\pi}\int_{x-1}^x\frac{\cos(\pi y)}{y\log^2(y)}dy\\
 & \sim & -\frac{2\cos(\pi x)}{\pi\log(x)}\quad\text{as } x\to\infty.
\end{eqnarray*}
Therefore, by Theorem~\ref{ginbv}(d), 
$f\ast g\in C^0(\Rbar)$ and $\lim_{|x|\to\infty}f\ast g(x)=0$
but for each $1\leq p<\infty$
we have $f\ast g\not\in L^p$.

\section{Differentiation and integration}
If $g$ is sufficiently smooth then the pointwise derivative is
$(f\ast g)'(x)=f\ast g'(x)$.  Recall the definition
$\acr$ of primitives of $L^1$ functions given in the proof of
Proposition~\ref{acr}.  In the following theorem we require
pointwise derivatives of $g$ to exist at each point in $\R$.
\begin{theorem}\label{difftheorem}
Let $f\in\alexc$, $n\in\N$ and $g^{(k)}\in AC(\Rbar)$ for each
$0\leq k\leq n$.  Then $f\ast g\in C^n(\R)$ and 
$(f\ast g)^{(n)}(x)=f\ast g^{(n)}(x)$ for each $x\in\R$.
\end{theorem}
\bigskip
\noindent
{\bf Proof:} First consider $n=1$.  Let $x\in\R$.  Then
\begin{equation}
(f\ast g)'(x)  =  \lim_{h\to 0}\intinf f(y)\left[\frac{g(x+h-y)
-g(x-y)}{h}\right]dy.\label{fg'}
\end{equation}
To take the limit inside the integral we can show that the bracketed
term in the integrand is of uniform bounded variation for
$0<|h|\leq 1$.  Let $h\not=0$.
Since $g\in\acr$ it  follows that the variation is given by the 
Lebesgue integrals 
\begin{align}
&V_{y\in\R}\left[\frac{g(x+h-y)
-g(x-y)}{h}\right]  =  \intinf \left|\frac{g'(x+h-y)
-g'(x-y)}{h}\right|dy\notag\\
 & \leq \intinf|g''(y)|dy+ \intinf \left|\frac{g'(x+h-y)
-g'(x-y)}{h}-g''(x-y)\right|dy.\label{L1der}
\end{align}
Since $g'\in\acr$  we have
$g''\in L^1$.  The  second integral on the right of \eqref{L1der}
gives the $L^1$ derivative of $g'$ in the limit $h\to 0$.
See \cite[p.~246]{folland}.  Hence, in \eqref{fg'} we
can use Theorem~22 in \cite{talviladenjoy} 
to take the limit under the integral sign.  This then
gives $(f\ast g)'(x)=f\ast g'(x)$.  Theorem~\ref{ginbv}(d) now
shows $(f\ast g)'\in C^0(\Rbar)$.  Induction on $n$
completes the proof.\qed

For similar results when $f\in L^1$,
see \cite[Proposition~8.10]{folland}.

Note that $g'\in\acr$ does not imply $g\in\acr$.  For example,
$g(x)=x$.  The conditions $g^{(k)}\in\bv$ for $0\leq k\leq n+1$
imply those in Theorem~\ref{difftheorem}. To see this it suffices
to consider $n=1$.  If $g', g''\in\bv$ then $g''$ exists at each
point and is bounded.  Hence, the
Lebesgue integral
$g'(x)=g'(0)+\int_0^x g''(y)\,dy$ exists for each $x\in\R$ 
and $g'$ is absolutely continuous.  Since $g'\in\bv$
we then have $g'\in\acr$.
Similarly for $n> 1$.
The example $g(x)=|x|^{1.5}\sin(1/[1+x^2])$
shows the $\acr$ condition
in the theorem is weaker that the above $\bv$ condition since
$g,g'\in\acr$ but $g''(0)$ does not exist so $g''\not\in\bv$.

We found that when $g\in\bv\cap L^1$ then $f\ast g\in\alexc$.
We can compute the the distributional derivative $(F\ast g)'=
f\ast g$ where $F$ is a primitive of $f$.
\begin{prop}\label{propdiff}
Let $F\in C^0(\Rbar)$ and write $f=F'\in\alexc$.  Let $g\in\bv\cap
L^1$.  Then $F\ast g\in C^0(\Rbar)$ and $(F\ast g)'=f\ast g\in\alexc$.
\end{prop}
\bigskip
\noindent
{\bf Proof:}
Let $x,t\in\R$.  Then by the usual H\"older inequality,
\begin{eqnarray*}
|F\ast g(x)-F\ast g(t)| & = & \left|\intinf \left[
F(x-y)-F(t-y)\right]g(y)\,dy\right|\\
 & \leq & \|F(x-\cdot)-F(t-\cdot)\|_\infty \|g\|_1\\
 & \to & 0\quad\text{as }t\to x \text{ since } F \text{ is 
uniformly continuous on } \R.
\end{eqnarray*}
Hence, $F\ast g$ is continuous on $\R$.  Dominated convergence shows
that $\lim_{x\to\pm\infty}F\ast g(x)=F(\pm\infty)\intinf g$.
Therefore, $F\ast g\in C^0(\Rbar)$.

Let $\phi\in\D$.  Then 
\begin{eqnarray*}
\langle (F\ast g)',\phi\rangle & = & -\langle F\ast g,\phi'\rangle\\
 & = & -\intinf\intinf F(x-y)g(y)\phi'(x)\,dy\,dx\\
 & = &  -\intinf g(y)\intinf F(x-y)\phi'(x)\,dx\,dy\text{
(Fubini--Tonelli theorem).} 
\end{eqnarray*}
Integrate by parts and use the change of variables $x\mapsto x+y$ to get
\begin{eqnarray*}
\langle (F\ast g)',\phi\rangle & = & \intinf g(y)\intinf f(x)\phi(x+y)\,dx
\,dy\\
 & = & \intinf f(x)\intinf g(y)\phi(x+y)\,dy\,dx \text{ (by Proposition~
\ref{fubini})}\\
 & = & \intinf f(x)\intinf g(y-x)\phi(y)\,dy\,dx\\
 & = & \intinf \phi(y)\intinf f(x)g(y-x)\,dx\,dy\text{ (by Proposition~
\ref{fubini})}\\
 & = & \langle f\ast g,\phi\rangle.\qed
\end{eqnarray*}

This gives an alternate definition of $f\ast g$ for $f\in\alexc$ and
$g\in L^1$.
\begin{theorem}
Let $f\in\alexc$, let $F\in\balexc$ be the
primitive of $f$ and let $g\in L^1$.
Define $f\ast g$ as in 
Definition~\ref{L1defn}.  Then $(F\ast g)'=f\ast g\in\alexc$.
\end{theorem}
\bigskip
\noindent
{\bf Proof:}
Let $-\infty<\alpha<\beta<\infty$.
Let $\{g_n\}\subset\bv\cap L^1$
such that $\|g_n-g\|_1\to 0$.  By Proposition~\ref{propdiff} we have
$$
\int_\alpha^\beta (F\ast g)'=F\ast g(\beta)-F\ast g(\alpha)
=\intinf F(y)\left[g(\beta-y)-g(\alpha-y)\right]\,dy.
$$  As
in \eqref{theorem1c},
$\int_\alpha^\beta f\ast g_n  = \intinf F(y)\left[
g_n(\beta-y)-g_n(\alpha-y)\right]\,dy$.  Hence,
\begin{align*}
&\left|\int_\alpha^\beta \left[(F\ast g)'-f\ast g_n\right]\right|\\
&=\left|\intinf F(y)\left[\left(g(\beta-y)-g_n(\beta-y)\right)
-\left(g(\alpha-y)-g_n(\alpha-y)\right)\right]\,dy\right|\\
&\leq \|F\|_\infty\left(\|g(\beta-\cdot)-g_n(\beta-\cdot)\|_1
+\|g(\alpha-\cdot)-g_n(\alpha-\cdot)\|_1\right)\\
&=2\|f\| \|g_n-g\|_1.
\end{align*}
Therefore, $\|(F\ast g)'-f\ast g_n\|\leq 
2\|f\| \|g_n-g\|_1 \to 0$  as $n\to\infty$. \qed

The next  theorem and its corollary give results on integrating convolutions.
\begin{theorem}
Let $f\in\alexc$ and let $g\in L^1$.  Define $F(x)=\int_{-\infty}^x f$
and $G(x)=\int_{-\infty}^x g$.  Then $f\ast G\in C^0(\Rbar)$ and
$f\ast G(x)=F\ast g(x)$ for all $x\in\R$.
\end{theorem}
\bigskip
\noindent
{\bf Proof:}
Since $G\in\acr$, Theorem~\ref{ginbv}(d) shows $f\ast G\in C^0(\Rbar)$.
We have 
\begin{eqnarray*}
f\ast G(x) & = & \intinf f(y)\int_{-\infty}^{x-y}g(z)\,dz\,dy\\
 & = & \intinf\intinf f(y)\chi_{(-\infty,x-y)}(z)g(z)\,dz\,dy\\
 & = & \intinf\intinf f(y)\chi_{(-\infty,x-y)}(z)g(z)\,dy\,dz\\
 & = & \intinf g(z)\int_{-\infty}^{x-z} f(y)\,dy\,dz\\
 & = & F\ast g(x).
\end{eqnarray*}
Proposition~\ref{fubini} justifies the interchange of orders of
integration.\qed
\begin{corollary}
(a) $f\ast g=(F\ast g)'=(f\ast G)'$
(b) For all $-\infty\leq \alpha<\beta\leq\infty$ we have
$\int_\alpha^\beta f\ast g=F\ast g(\beta)-F\ast g(\alpha)
=f\ast G(\beta)-f\ast G(\beta)$.
\end{corollary}
Hence, the convolution $f\ast g$ can be evaluated by taking
the distributional derivative of the Lebesgue integral $F\ast g$.
Since $f\ast G\in C^0(\Rbar)$ when $f\in\alexc$ and $G\in \bv$
we can use the equation $f\ast g=(f\ast G)'$ to define $f\ast g$
for $f\in\alexc$ and $g=G'$ for $G\in\bv$.  In this case, $g$ will be
a signed Radon measure.  As $G(x)=\int_{-\infty}^x g$ and this integral
is a regulated primitive integral \cite{talvilarpi}, we will save
this case for discussion elsewhere.

\section{Appendix}
The integration  by parts formula is as follows.  If $f\in\alexc$
and $g\in\bv$  it  gives the integral of $fg$ in terms  of 
a Henstock--Stieltjes integral:
\begin{eqnarray}
\int_{-\infty}^\infty fg & = &
F(\infty)g(\infty)-\int_{-\infty}^\infty 
F\,dg.\label{parts}
\end{eqnarray}
See \cite{talviladenjoy} and \cite[p.~199]{mcleod}.

We have the following corollary for functions of essential
bounded variation.
\begin{corollary}\label{partscorollary}
Let $F\in C^0(\Rbar)$.  Let $g\in\ebv$.  Fix $0\leq \gamma\leq 1$.
Take $g_\gamma\in\nbv_\gamma$ such that $g_\gamma=g$ almost
everywhere.  Let $\mu_g$ be the signed Radon measure given by
$g'$.  Then $\intinf F\,dg_\gamma=\intinf F\,d\mu_g$.
\end{corollary}
\bigskip
\noindent
{\bf Proof:}
The distributional derivative of $g$ is
$\langle g',\phi\rangle=-\langle g,\phi'\rangle=-\intinf g\phi'=
\intinf \phi\,d\mu_g$ for all $\phi\in\D$.  Note that $g_\gamma$
is unique and $\mu_g=\mu_{g_\gamma}$.  Suppose $\phi\in\D$ with
${\rm supp}(\phi)\subset[A,B]\subset\R$.  Then, using integration
by parts for the Henstock--Stieltjes integral,
\begin{eqnarray*}
\langle g_\gamma,\phi'\rangle & = & \int_A^Bg_\gamma\phi'
 =  g_\gamma(B)\phi(B)-g_\gamma(A)\phi(A)-\int_A^B\phi\,
dg_\gamma
  =  -\intinf\phi\,
dg_\gamma\\
 & = & -\langle g_\gamma',\phi\rangle=-\intinf\phi\,d\mu_{g_\gamma}
  =  -\intinf\phi\,d\mu_g.
\end{eqnarray*}

Let $F\in C^0(\Rbar)$.  There is a uniformly bounded sequence
$\{\phi_n\}\subset\D$ such that $\phi_n\to F$ pointwise on $\R$.
By dominated convergence,
$$
\lim_{n\to\infty}\intinf\phi_n\,dg_\gamma=\intinf F\,dg_\gamma
=\lim_{n\to\infty}\intinf\phi_n\,d\mu_g=\intinf F\,d\mu_g.\qed
$$

Corollary~\ref{partscorollary} now justifies the following definition.
\begin{defn}\label{defnebv}
Let $f\in\alexc$ and let $F\in\balexc$ be its primitive. Let $g\in\ebv$.  
Fix $0\leq\gamma\leq 1$ and
take $g_\gamma\in\nbv_\gamma$ such that $g_\gamma=g$
almost everywhere.  Define
$$
\intinf fg = g_\gamma(\infty)F(\infty)-\intinf F\,d\mu_g =\intinf f
g_\gamma.
$$
\end{defn}
Since limits at infinity are not affected by the choice of
$\gamma$, the definition is independent of $\gamma$.

The H\"older inequality is
\begin{equation}
\left|\int_{-\infty}^\infty
fg\right|\leq \left|\intinf f\right|\inf_{\R}|g|+\|f\|Vg\leq
\|f\| \|g\|_{\bv}\label{holder}
\end{equation}
and is valid for all $f\in\alexc$ and $g\in\bv$.  For
$g\in\ebv$ we replace $g$ with $g_\gamma$.  This
gives 
\begin{equation}
\left|\int_{-\infty}^\infty
fg\right|\leq \left|\intinf f\right|\inf_{\R}|g_\gamma|+\|f\|Vg_\gamma\leq
\|f\| \|g\|_{\ebv}.\label{ebvholder}
\end{equation}
See 
\cite[Lemma~24]{talvilafourier}
for a proof using the Henstock--Kurzweil integral.
The same proof  works for the continuous primitive integral.

A Fubini theorem has been established in
\cite{ang} for the continuous primitive integral on  compact intervals.  
This says that
if a double integral exists in the plane then the two iterated integrals exist
and are equal.  Of more utility for the case at hand is to show directly
that iterated integrals are equal  without resorting to the double integral.
The following theorem extends a type of Fubini theorem  proved  on page~$58$ in
\cite{celidze} for
the wide Denjoy integral on compact intervals.
\begin{prop}\label{fubini}
Let $f\in\alexc$.
Let $g\fn\R\times\R\to\R$ be measurable.
Assume (i) for each $x\in\R$ the function
$y\mapsto g(x,y)$  is in $\bv$; (ii) the function $x\mapsto
V_{y\in\R}g(x,y)$ is in $L^1$; (iii) there is $M\in L^1$
such that for each $y\in \R$ we have $|g(x,y)|\leq M(x)$.  Then
the iterated integrals exist and are equal,
$\intinf\intinf f(y)g(x,y)\,dy\,dx=
\intinf\intinf f(y)g(x,y)\,dx\,dy$.
\end{prop}
\bigskip
\noindent
{\bf Proof:} Let $F\in\balexc$ be the primitive of $f$.  For
each $x\in\R$ we have
$$
\intinf f(y)g(x,y)\,dy  =  F(\infty)g(x,\infty)-\intinf F(y)\,d_2g(x,y)
$$
where $d_2(x,y)$ indicates a Henstock--Stieltjes integral with respect
to $y$.  Then,
\begin{eqnarray}
\intinf\intinf f(y)g(x,y)\,dy\,dx & = & F(\infty)\intinf g(x,\infty)\,dx\notag\\
 & & \quad-
\intinf\intinf F(y)\,d_2g(x,y)\,dx.\label{fubiniproof1}
\end{eqnarray}
The integral $\intinf g(x,\infty)\,dx$ exists due
to condition (iii).
The iterated  integral in \eqref{fubiniproof1}
converges absolutely since
\begin{eqnarray*}
\left|\intinf\intinf F(y)\,d_2g(x,y)\,dx\right| & \leq &
\intinf\left|\intinf F(y)\,d_2g(x,y)\right|dx\\
 & \leq & \|F\|_\infty\intinf V_{y\in\R}g(x,y)\,dx.
\end{eqnarray*}

Now, show the function $y\mapsto \intinf g(x,y)\,dx$ is in $\bv$.
Let $\{(s_i,t_i)\}_{i=1}^n$ be disjoint intervals in $\R$.  Then
\begin{eqnarray*}
\sum_{i=1}^n\left|\intinf g(x,s_i)\,dx - \intinf g(x,t_i)\,dx\right| 
& \leq & 
\sum_{i=1}^n\intinf \left|g(x,s_i)-g(x,t_i)\right| dx\\
 & = & \intinf\sum_{i=1}^n \left|g(x,s_i)-g(x,t_i)\right| dx\\
 & \leq & \intinf V_{y\in\R}g(x,y)\,dx.
\end{eqnarray*}
The interchange of summation  and integration follows from condition
(ii) and the usual Fubini--Tonelli theorem.  Hence, the function
$y\mapsto \intinf g(x,y)\,dx$ is in $\bv$ and the iterated integral
$\intinf f(y)\intinf g(x,y)\,dx\,dy$ exists.

Integrate by parts,
\begin{eqnarray}
\intinf f(y)\intinf g(x,y)\,dx\,dy & = & F(\infty)\intinf g(x,\infty)\,dx
\label{fubiniproof2}\\
 & & \qquad
-\intinf F(y)\,d\left[\intinf g(x,y)\,dx\right].
\label{fubiniproof3}
\end{eqnarray}
In \eqref{fubiniproof2}, we have 
$\lim_{y\to\infty}\intinf g(x,y)\,dx=\intinf g(x,\infty)\,dx$
due to dominated convergence and condition (iii).  To complete the
proof we need to show the integrals in \eqref{fubiniproof1}
and \eqref{fubiniproof3} are equal.  First consider the case when
$F=\chi_{(a,b)}$ for an interval $(a,b)\subset\R$.  Then \eqref{fubiniproof1} becomes
$
\intinf\int_a^b d_2g(x,y)\,dx = \intinf\left[ g(x,b) -g(x,a)\right]dx
$.
And now \eqref{fubiniproof3} becomes
$
\int_a^b d\left[\intinf g(x,y)\,dx\right] =
\intinf g(x,b)\,dx -\intinf g(x,a)\,dx
$.
Hence, when $F$ is a step function, $F(y)=\sum_{i=1}^n c_i\chi_{I_i}(y)$
for some $n\in\N$, 
disjoint intervals $\{I_i\}_{i=1}^n$ and real numbers $\{c_i\}_{i=1}^n$,
we have the desired equality of \eqref{fubiniproof1} and \eqref{fubiniproof3}.
But $F\in\balexc$ is uniformly continuous on $\Rbar$, i.e., for each
$\epsilon>0$ there is $\delta >0$ such that for all $0\leq|x-y|<\delta$
we have $|F(x)-F(y)|<\epsilon$, for all $x<-1/\delta$  we have
$|F(x)|<\epsilon$ and for  all $x>1/\delta$ we have $|F(x)-F(\infty)|<
\epsilon$.  It then  follows from the compactness of  $\Rbar$ that
the step functions are dense in $\balexc$.  Hence, there is a sequence
of step functions $\{\sigma_N\}$ such that $\|F-\sigma_N\|_\infty\to 0$.
In \eqref{fubiniproof1} we have 
$$
\lim_{N\to\infty}\intinf\intinf \sigma_N(y)\,d_2g(x,y)\,dx
  =  \intinf\intinf F(y)\,d_2g(x,y)\,dx.
$$ The $N$ limit  can be brought inside the $x$ integral using dominated
convergence and  (ii) since
$|\intinf \sigma_N(y)\,d_2g(x,y)|\leq (\|F\|_\infty+1)V_{y\in\R}g(x,y)$
for large enough $N$.  The $N$ limit  can be brought inside the $y$ 
integral  using dominated convergence  since 
$|\sigma_N(y)|\leq (\|F\|_\infty+1)$ for large enough $N$.
In \eqref{fubiniproof3} we have
$$
\lim_{N\to\infty}\intinf \sigma_N(y)\,d\left[\intinf g(x,y)\,dx\right] 
=
\intinf F(y)\,d\left[\intinf g(x,y)\,dx\right].
$$
The $N$ limit can be brought inside the $y$ integral since
$\{\sigma_N\}$ converges to $F$ uniformly on $\Rbar$
and $d\left[\intinf g(x,y)\,dx\right]$ is a finite signed measure.
\qed
\begin{corollary}\label{corfubini}
If $f$ has compact support we can replace (iii) with:
(iv) for each $y\in {\rm supp}(f)$ the function $x\mapsto g(x,y)$
is in $L^1$.
\end{corollary}

\end{document}